\documentclass[10pt,journal]{IEEEtran}

\usepackage[utf8]{inputenc} % allow utf-8 input
\usepackage[T1]{fontenc}    % use 8-bit T1 fonts
\usepackage{hyperref}       % hyperlinks
\usepackage{url}            % simple URL typesetting
\usepackage{booktabs}       % professional-quality tables
\usepackage{amsfonts}       % blackboard math symbols
\usepackage{nicefrac}       % compact symbols for 1/2, etc.
\usepackage{microtype}      % microtypography
\usepackage{amsmath,cite,xspace}
\usepackage{graphicx}
\usepackage{epsfig}
\usepackage{epstopdf}
\usepackage{makecell}
\usepackage{bm}
\newcommand{\subparagraph}{}
\usepackage{color}
\usepackage{amssymb}
\usepackage{hyperref}
\newcommand{\utwi}[1]{\mbox{\boldmath $ #1$}}
\usepackage[table,xcdraw]{xcolor}


\usepackage{mathtools}

\newcommand{\reals}{\mathbb{R}}

\newcommand{\bx}{{\bf x}}

\newcommand{\bv}{{\bf v}}

\newcommand{\bP}{{\bf P}}

\newcommand{\bQ}{{\bf Q}}

\newcommand{\btheta}{{\utwi{\theta}}}

\graphicspath{{./Figures/}} \DeclareGraphicsExtensions{.eps,.pdf,.jpg,.mps,.png}

\title{Emulating AC OPF solvers for Obtaining Sub-second Feasible, Near-Optimal Solutions}

 \author{\IEEEauthorblockN{Kyri Baker}\\
 \IEEEauthorblockA{University of Colorado Boulder\\
Renewable and Sustainable Energy Institute\\
kyri.baker@colorado.edu}
}
\makeatletter
\let\old@ps@headings\ps@headings
\let\old@ps@IEEEtitlepagestyle\ps@IEEEtitlepagestyle
\def\psccfooter#1{%
    \def\ps@headings{%
        \old@ps@headings%
        \def\@oddfoot{\strut\hfill#1\hfill\strut}%
        \def\@evenfoot{\strut\hfill#1\hfill\strut}%
    }%
    \def\ps@IEEEtitlepagestyle{%
        \old@ps@IEEEtitlepagestyle%
        \def\@oddfoot{\strut\hfill#1\hfill\strut}%
        \def\@evenfoot{\strut\hfill#1\hfill\strut}%
    }%
    \ps@headings%
}
\makeatother

\begin{document}
\maketitle

\begin{abstract}
Using machine learning to obtain solutions to AC optimal power flow has recently been a very active area of research due to the astounding speedups that result from bypassing traditional optimization techniques. However, generally ensuring feasibility of the resulting predictions while maintaining these speedups is a challenging, unsolved problem. In this paper, we train a neural network to emulate an iterative solver in order to cheaply and approximately iterate towards the optimum. Once we are close to convergence, we then solve a power flow to obtain an overall AC-feasible solution. Results shown for networks up to 1,354 buses indicate the proposed method is capable of finding feasible, near-optimal solutions to AC OPF in milliseconds on a laptop computer. In addition, it is shown that the proposed method can find ``difficult'' AC OPF solutions that cause flat-start or DC-warm started algorithms to diverge.
\end{abstract}

\section{Introduction}
AC optimal power flow (OPF) is a canonical power systems operation problem that is at the heart of optimizing large-scale power networks. Solving this problem quickly and efficiently has been the subject of decades of research. One particularly interesting development in this area is the use of machine learning (ML), in particular deep learning, to obtain solutions to AC OPF \cite{Zamzam_learn_19,chatzos2020highfidelity,DeepOPF_AC}. 
Within this area, ensuring feasibility of the resulting solution has been a challenge. In this paper, we propose a neural network model which aims to emulate an AC OPF solver (in particular, the Matpower Interior Point Solver (MIPS), although the framework is not specific to MIPS). The benefit of using an ML model instead of the MIPS solver directly is that no matrix inverses or factorizations are needed, and inference is extremely fast, resulting in an overall faster convergence. While we do not claim that feasibility can be guaranteed for every single output of the learning-based model, empirically, we have observed very positive results on the chosen networks in terms of optimality gap, speed, and convergence success. 
\begin{figure}[t!]
    \centering
    \hspace*{-2mm}\includegraphics[width=0.52\textwidth]{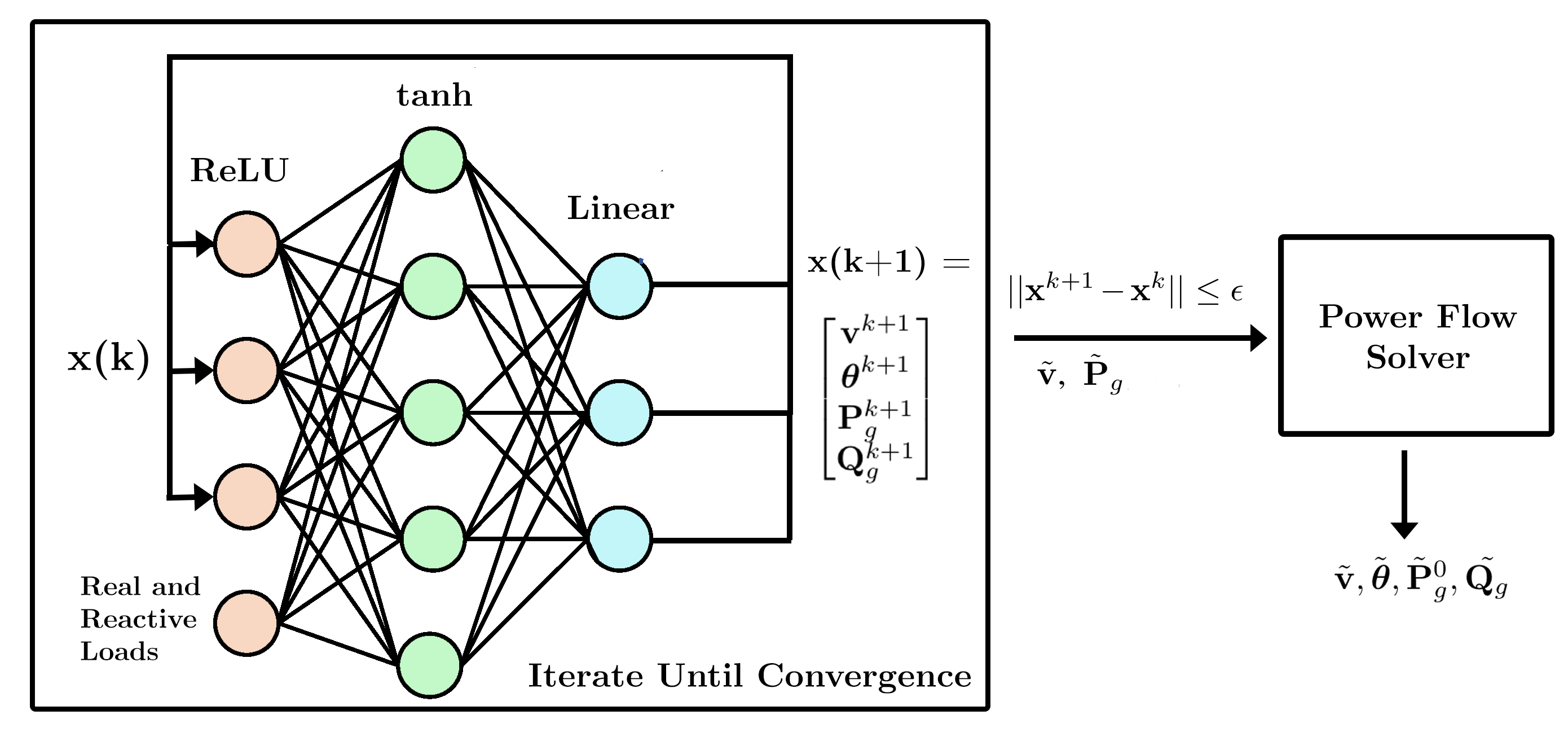}
    \caption{Using a neural network to approximate fast iterations towards the optimum, then solving a power flow to recover feasibility.}
    \label{fig:RNN}\vspace{-12pt}
\end{figure}
The model proposed here is comprised of a fully connected three-layer neural network (NN) $F_R$ with feedback, where input $\bx^{k}$ is the candidate optimal solution vector at iteration $k$. Reminiscent of a simple recurrent neural network, the model iteratively uses feedback from the output layer as inputs until convergence ($||\bx^{k+1}-\bx^{k}|| \leq \epsilon$). The model thus bypasses any construction of a Jacobian matrix or associated inverse, for example. This has benefits when trying to solve challenging AC OPF problems that may encounter near-singular power flow Jacobians while solving. In order to obtain a final AC feasible solution, a subset of the learned variables are sent to a power flow solver. See Fig. \ref{fig:RNN} for an overview of the testing phase of the algorithm, where the tilde over the variables indicates the candidate optimal solution produced by the model.

Some works have looked at penalizing constraint violations during training, which can help preserve AC feasibility, but cannot guarantee it \cite{chatzos2020highfidelity, DeepOPF_AC}. In addition, these techniques can result in very cumbersome-to-train loss functions. Other work using ML for AC OPF has recovered AC feasibility by using a post-processing procedure with the AC power flow equations, but requires a restricted training set generated from a modified AC OPF problem, sometimes requires PV/PQ switching, and was only tested on small networks \cite{Zamzam_learn_19}. The approach presented in this paper utilizes a similar concept to the latter technique, but offers advantages over all previous techniques to ensure feasibility:

\begin{itemize}
    \item No restriction on the training set is required; previous AC OPF runs can be used to train the ML model.
    \item The model \emph{emulates} an iterative algorithm, meaning that each model run is a small step towards the optimum, instead of directly predicting the AC OPF solution.
    \item While only a subset of variables is sent to the power flow solver, the ML model utilizes information about the entire OPF solution, better informing the model as it iterates towards the optimum. Evolutions of variable iterations towards the optimum can also be observed while using the model, unlike completely black-box approaches.
    \item A deep network is not required for this framework; in this paper, a shallow neural network (only one hidden layer) is used, which facilitates very fast training (under 12 hours on a laptop) and simpler parameter tuning.
\end{itemize}

These prior works \cite{Zamzam_learn_19,chatzos2020highfidelity,DeepOPF_AC} additionally only perturb the given loads by 10\% or 20\% in both training and testing, which results in a much smaller set of optimal solutions, making it easier for the ML model to map these conditions onto optimal values. Here, we generate data than spans a larger region of the feasible space, including solutions which represent the state of the system near voltage collapse. These situations pose challenging situations for traditional AC OPF solvers, which will be discussed in more detail.

Promising initial results are shown for 30, 500, and 1,354-bus networks suggesting that the learning-based model may provide a desirable tradeoff between speed and optimality gap while maintaining feasibility.

% \section{Optimal Power Flow}
% To provide the necessary background to discuss the feasibility of the DC OPF solutions in the AC OPF problem, we first briefly summarize three of the most common ways of optimizing generation dispatch in transmission networks. First, define coefficients $a_j$, $b_j$, and $c_j$ as the operational costs associated with generator $j$. Let set $\cal{G}$ be the set of all generators in the network, $\cal{N}$ be the set of all nodes (buses) in the network, $\cal{L}$ be the set of all lines (branches) in the network, $\Omega_i$ be the set of buses connected to bus $i$, and $\cal{G}_i$ be the set of generators connected to bus $i$. Define $p_{l,j}$ ($q_{l,j}$) as the total active (reactive) power consumption at node $j$, $p_{g,j}$ ($q_{g,j}$) is the active (reactive) power output of generator $j$, and $\underline{p}_{g,j}$ ($\underline{q}_{g,j}$) and $\overline{p}_{g,j}$ ($\overline{q}_{g,j}$) are lower and upper limits on active (reactive) power generation, respectively. The complex voltage at node $i$ has magnitude $|v_i|$ and phase angle $\theta_i$, and the difference in phase angle between neighboring buses $i$ and $m$ is written as $\theta_m$. Constant parameters $G_{im}$ and $B_{im}$ are the conductance and susceptance of line $im$, respectively.

\section{Iterating using inference} \label{sec:LBQN}
A general nonconvex optimization problem with $n$-dimensional optimization variable vector $\bx$, cost function $f(\cdot): \reals^{n} \rightarrow \reals$, $M$ equality constraints $g_i(x) = 0$, $g_i(\cdot): \reals^{n} \rightarrow \reals$, and $P$ inequality constraints $h_j(x) \leq 0$, $h_j(\cdot): \reals^{n} \rightarrow \reals$ can be written as

\begin{subequations} \label{eqn:gen_opt}
\begin{align} 
\min_{\substack{\bx}} ~~ f(\bx)& \\
\mathrm{s.t:~} g_i(\bx) &= 0, ~~i=1, ... M\\
h_j(\bx) &\leq 0, ~~j=1, ..., P
 \end{align}
\end{subequations}

Many iterative optimization solvers use Hessians of the Lagrangian function to iterate towards the optimal solution of constrained nonconvex problems, including the Matpower Interior Point Solver (MIPS) \cite{MATPOWER}, which leverages a primal-dual interior point algorithm to update candidate solution $\bx^{k}$ at iteration $k$. Instead of using Lagrangian functions or forming Hessian matrices, the learning-based method uses a deep learning model $F_R(\cdot): \reals^n \rightarrow \reals^n$ that takes in $\bx^{k}$ as an input and provides $\bx^{k+1}$ as an output; e.g.

\begin{equation}\label{eqn:LBNR}
    \bx^{k+1} = F_R(\bx^{k}).
\end{equation}

A fully connected three-layer NN is used here. The variable vector $\bx^{k} = [\bv^{k},~ \btheta^k,~ \bP^k_g,~ \bQ^k_g]^T$, where the iteration index is $k$, $\bv$ contains the complex voltage magnitudes at each bus, $\btheta$ contains the complex voltage angles, and $\bP_g$ and $\bQ_g$ are the real and reactive power outputs at each generator, respectively. 

\subsection{Network architecture}
A rectified linear unit (ReLU) was used as the activation function on the input layer; a hyperbolic tangent (tanh) activation function was used in the hidden layer, and a linear function was used on the output. Bounds on generation and voltages were enforced with a threshold on the output layer of the NN. Normalization of the inputs was also performed such that all inputs were in the range $[0,1]$, which improved NN performance. In addition to $\bx^{k}$, the network loading (constant throughout inference) was given as an input to the NN. The chosen number of nodes in the hidden layer and training dataset sizes are shown in Table \ref{tab:hyper}. The data generation, training and testing of the network, and simulations were performed on a 2017 MacBook Pro laptop with 16 GB of memory. Keras with the Tensorflow backend was used to train the neural network using the Adam optimizer.

\begin{table}[t!]\centering\small \caption{Number of nodes and training samples}
\begin{tabular}{|l|l|l|l|l|}
\hline
\textbf{Case}      & \textbf{\begin{tabular}[c]{@{}l@{}}Input\\ Nodes\end{tabular}} & \textbf{\begin{tabular}[c]{@{}l@{}}Output\\ Nodes\end{tabular}} & \textbf{\begin{tabular}[c]{@{}l@{}}Hidden\\ Nodes\end{tabular}} & \textbf{\begin{tabular}[c]{@{}l@{}}Training\\ Samples\end{tabular}} \\ \hline
\textbf{30-bus}    & 112            & 72             & 100   & 72,111    \\ \hline
\textbf{500-bus} & 1,512          & 1,112          & 2,300 & 111,674 \\ \hline
\textbf{1,354-bus} & 4,470          & 3,228         & 6,000 & 126,724 \\ \hline
\end{tabular}\label{tab:hyper}
\end{table}

\subsection{Data Generation}
The MATPOWER Interior Point Solver (MIPS) \cite{MATPOWER} was used to generate the data and was used as the baseline for comparison with the NN model. A single training sample consists of the pair $[\bx^{k}, \bx^{k+1}]$ obtained from the solver. The termination tolerance of the MIPS solver was set to $10^{-9}$ for data generation and $10^{-4}$ for testing. The tolerance of the learning-based solver was set to $10^{-4}$, where convergence is reached when $||\bx^{k+1}-\bx^k|| \leq \epsilon$. A smaller tolerance was used for data generation to promote smoother convergence and ``basins of attraction" within the ML model. For a fair comparison, the same convergence criteria was used for the NN model and MIPS during testing. 500 different loading scenarios were randomly generated at each load bus from a uniform distribution of $\pm40\%$ around the given base loading scenario in MATPOWER. Table \ref{tab:hyper} shows the number of training samples generated for each scenario. Each generated set of loads was used to solve a standard power flow, and if the power flow could not find a solution, the data was not included in the training set. It is recognized that generating a diverse and representative dataset is an important and essential thrust of research within learning-based OPF methods. This is an important direction of future work.

\begin{table}[t!]\centering \caption{The ML solution finds feasible AC OPF solutions faster than traditional methods and solves in a comparable time to DC OPF.}
\begin{tabular}{|l|l|l|l|}
\hline
\textbf{Network}   & \textbf{OPF Type} & \textbf{\begin{tabular}[c]{@{}l@{}}Average Solve\\ Time (s)\end{tabular}} & \textbf{\begin{tabular}[c]{@{}l@{}}Percent Solve\\ Success (\%)\end{tabular}} \\ \Xhline{2\arrayrulewidth}
\textbf{}          & \textbf{\begin{tabular}[c]{@{}l@{}}AC OPF \\ w/Flat Start\end{tabular}}   & 0.069 s   & 100\%  \\ \cline{2-4} 
 & \textbf{\begin{tabular}[c]{@{}l@{}}AC OPF \\ w/DC Start\end{tabular}}        & 0.079 s    & 100\%     \\ \cline{2-4} 
\textbf{30-bus}    & \textbf{\begin{tabular}[c]{@{}l@{}}AC OPF \\ w/PF Start\end{tabular}}     & 0.068 s & 100\%     \\ \cline{2-4}   & \cellcolor[HTML]{EFEFEF}\textbf{\begin{tabular}[c]{@{}l@{}}NN with PF\end{tabular}} & \cellcolor[HTML]{EFEFEF}0.050 s   & \cellcolor[HTML]{EFEFEF}100\%  \\ \cline{2-4}   & \textbf{DC OPF}   & 0.010 s & 100\%  \\ \cline{2-4}   & \textbf{DC OPF w/PF}   & 0.018 s & 100\%  \\\Xhline{2\arrayrulewidth}
\textbf{}          & \textbf{\begin{tabular}[c]{@{}l@{}}AC OPF\\ w/Flat Start\end{tabular}}   & 3.26 s & 100\%     \\ \cline{2-4}   & \textbf{\begin{tabular}[c]{@{}l@{}}AC OPF \\ w/DC Start\end{tabular}}         & 3.34 s        & 100\%   \\ \cline{2-4} 
\textbf{500-bus}   & \textbf{\begin{tabular}[c]{@{}l@{}}AC OPF \\ w/PF Start\end{tabular}}   & 3.15 s& 100\%   \\ \cline{2-4}  & \cellcolor[HTML]{EFEFEF}\textbf{\begin{tabular}[c]{@{}l@{}}NN with PF \end{tabular}} & \cellcolor[HTML]{EFEFEF}0.12 s    & \cellcolor[HTML]{EFEFEF}100\%     \\ \cline{2-4} 
    & \textbf{DC OPF}  & 0.12 s  & 100\%  \\ \cline{2-4} 
    & \textbf{DC OPF w/PF}  & 0.13 s  & 100\% \\\Xhline{2\arrayrulewidth}
\textbf{}          & \textbf{\begin{tabular}[c]{@{}l@{}}AC OPF\\ w/Flat Start\end{tabular}} & 20.12 s  & 59.4\%    \\ \cline{2-4}   & \textbf{\begin{tabular}[c]{@{}l@{}}AC OPF \\ w/DC Start\end{tabular}} & 11.63 s& 59.4\%    \\ \cline{2-4} 
\textbf{1354-bus} & \textbf{\begin{tabular}[c]{@{}l@{}}AC OPF \\ w/PF Start\end{tabular}}   & 10.64 s & 59.8\%    \\ \cline{2-4}  & \cellcolor[HTML]{EFEFEF}\textbf{\begin{tabular}[c]{@{}l@{}}NN with PF\end{tabular}} & \cellcolor[HTML]{EFEFEF}0.25 s      & \cellcolor[HTML]{EFEFEF}98.00\%      \\ \cline{2-4}  & \textbf{DC OPF}   & 0.33 s  & 100\%    \\ \cline{2-4}  & \textbf{DC OPF w/PF}   & 0.38 s  & 99.20\%       \\ \hline
\end{tabular}\label{tab:time}
\end{table}

\section{Simulation Results} \label{sec:sim_result}
Here, three networks were considered: The IEEE 30 bus, IEEE 500 bus, and 1,354-bus PEGASE networks \cite{PEGASE}. Both training and testing were performed locally on a MacBook Pro with 16 GB of RAM. To be consistent with the constraint set across networks, because some of the considered networks do not contain line flow limits, these were neglected in all networks. However, it is not expected that the inclusion of these constraints would dramatically change the results shown here. In Table \ref{tab:time}, the learning-based method (``NN'') was compared with 4 other cases across the 500-sample training set: AC OPF flat-started, AC OPF warm-started with both a DC OPF solution and a power flow solution, a DC OPF, and a DC OPF with an AC PF in post-processing to pursue AC feasibility. While the DC OPF never produces an AC feasible solution, it is often used as an approximation for AC OPF and used in many LMP-based markets to calculate prices and thus provides an interesting comparison. Note that in these simulations, no ``extra'' feasibility steps were performed (like PV/PQ switching in \cite{Zamzam_learn_19}); the output of the NN was directly sent to a power flow solver. 

\subsection{Speedups} As Table \ref{tab:time} shows, especially for the 500 and 1,354-bus networks, the learning-based solution can provide solutions with speedups of over twenty times faster than using a traditional solver, even one that has been warm-started. A DC OPF, in comparison, takes about the same amount of time to solve a much simpler, convex problem, and, considering it entails a linear set of constraints, has no issues finding a solution. For smaller networks, the benefit of using a ML-based method is negligible; sub-second solution times are already achieved by standard AC OPF methods. Another interesting idea is to use a ML model to warm-start the AC OPF as in \cite{Baker_learning} or the AC PF as in \cite{chen2020hotstarting}, but we only compare with traditional warm-start methods here. Lastly, it should be mentioned that the NN-based solver can be trained and tested using data generated from commercial solvers such as CPLEX or Gurobi which may be able to improve the performance of all of the test cases here. As the training data for the NN is simply the variable vector at each iteration, any iterative algorithm can be used for data generation. We expect the \emph{relative} speedups to remain relatively consistent regardless of the solver, however.

\subsection{Challenging OPF scenarios} In the 1,354-bus test case, load profiles were generated that challenged the MIPS solver. Typically for these difficult cases, continuation methods can be used to robustly solve the AC OPF by solving a series of simpler OPF problems \cite{homoSanja}. While robust, these methods can be very time consuming and may not be suitable for real-time operation. The learning-based method may provide an alternative to these methods, as issues with singular Jacobian matrices and solutions close to voltage instability do not affect inference as much. These solutions may still prove challenging for the post-processing power flow step, however, which is perhaps why a few failures were still encountered when using the learning-based method. 

In Table \ref{tab:time}, the last column refers to the percentage of trial runs in which the solver converged to a feasible solution. For the AC OPF cases, this refers to successful convergence of the MIPS AC OPF solver to point which satisfies the constraints. For the NN and DC OPF w/PF case, this refers to the MIPS AC PF solver converging to a point which satisfies the AC power flow constraints. Lastly, for the DC OPF case, this refers to the percentage of runs in which the MIPS solver found the optimal solution to the DC OPF problem. Considering the DC OPF problem is convex (and linear in some cases where the generator costs are linear), if the initial loading point was feasible, it is reasonable that the percent solve success is high.

\begin{table}[t!]\centering \small \caption{Percent error for the proposed method and for DC OPF across the 500-sample testing sets.}
\begin{tabular}{|l|l|l|l|l|}
\hline
\textbf{Network}   & \textbf{\begin{tabular}[c]{@{}l@{}}Average  \\\% Error \\(NN)\end{tabular}} & \textbf{\begin{tabular}[c]{@{}l@{}}Average \\\% Error \\(DC)\end{tabular}} & \textbf{\begin{tabular}[c]{@{}l@{}}Worst \\ \% Error \\(NN)\end{tabular}} & \textbf{\begin{tabular}[c]{@{}l@{}}Worst  \\ \% Error \\(DC)\end{tabular}} \\ \hline
\textbf{30-bus}    & 0.20\%     & 1.46\%       & 0.31\%             & 1.82\%    \\ \hline
\textbf{500-bus}   & 2.70\%  & 5.30\%   & 10.12\%       & 16.22\%   \\ \hline \textbf{1354-bus} & 0.09\%                  & 1.40\%     & 0.20\%                    & 1.73\%     \\ \hline
\end{tabular}\label{tab:optgap}
\end{table}

\begin{table}[t!]\centering \small \caption{Percent error for the proposed method and for DC OPF with an AC PF across the 500-sample testing sets.}
\begin{tabular}{|l|l|l|l|l|}
\hline
\textbf{Network}   & \textbf{\begin{tabular}[c]{@{}l@{}}Average  \\\% Error \\(NN)\end{tabular}} & \textbf{\begin{tabular}[c]{@{}l@{}}Average \\\% Error \\(DC w/PF)\end{tabular}} & \textbf{\begin{tabular}[c]{@{}l@{}}Worst \\ \% Error \\(NN)\end{tabular}} & \textbf{\begin{tabular}[c]{@{}l@{}}Worst  \\ \% Error \\(DC w/PF)\end{tabular}} \\ \hline
\textbf{30-bus}    & 0.20\%     & 0.26\%       & 0.31\%             & 0.37\%    \\ \hline
\textbf{500-bus}   & 2.70\%  & 29.13\%   & 10.12\%       & 35.34\%   \\ \hline \textbf{1354-bus} & 0.09\%                  & 1.87\%     & 0.20\%                    & 2.22\%     \\ \hline
\end{tabular}\label{tab:optgap2}
\end{table}

\subsection{Error of the learning-based solution vs. MIPS} Tables \ref{tab:optgap} and \ref{tab:optgap2} show the difference in cost for the learning-based solution (``NN'') and DC OPF without and with an AC power flow (respectively), calculated as follows:

\begin{align}\label{eqn:error}
\%~ \text{Error} = \Big|\frac{f(\bx^*) - f(\Tilde{\bx})}{f(\bx^*)}\Big|\cdot 100\%.
\end{align}

\noindent Where $f(\cdot)$ is the AC OPF objective function, and $\bx^*$ is the solution from MIPS, and $\Tilde{\bx}$ is the solution from the method being compared (i.e., the learning based or the DC solution). The equation in \eqref{eqn:error} is akin to the definition of ``optimality gap''; however, we avoid using this specific term here because there is no guarantee that $f(\bx^*)$ is the globally optimal solution. Due to the nonconvex nature of the problem and the fact that we are training and testing the model using a solver that is not a global optimization solver, assessing the true optimality gap is a direction of future work.

In the results, the DC OPF solution with and without AC PF post processing is compared against the AC OPF solution provided by MIPS. This is a key comparison because in many markets, DC OPF is used to calculate prices; thus, it is currently deemed an acceptable way of approximating the cost of network operation. As the table shows, however, the learning-based method produces even lower operating costs than DC OPF, on average and in the worst case throughout the testing set. Surprisingly, the solution obtained with the AC PF post-processing step, while now AC feasible, were found to overall increase the cost difference with the AC OPF solution. This may be due to the fact that when a power flow is run, all generator values except the slack bus generator are fixed to the given values (in this case, the DC OPF solution). Due to the fact that the DC OPF does not account for losses, this means that the slack generator must account for \emph{all} additional power requirements to obtain AC feasibility. This generator may not be the cheapest / best generator to serve this purpose; thus, the objective value ends up increasing with this additional power generation at the cost of feasibility.

\section{Comparison with other power flow models} \label{sec:comparison}
Even though DC OPF is used in many existing network operations, it is useful to compare the proposed method against other convexifications. In this section, the learning-based method is compared against the popular second-order cone programming relaxation (SOCP) in terms of speed and accuracy. The particular implementation that will be used for the comparison is that found in the \texttt{PowerModels} package \cite{powermodels} for SOCP \cite{SOCPJabr}. Table \ref{tab:socp} shows a side-by-side comparison of solving SOCP using IPOPT in Julia compared to the learning-based method, across the same testing set. While the learning based method appears to be faster and more accurate, it is important to note that unlike SOCP, it does not provide a theoretical lower bound for the true optimum. Additionally, an AC power flow could be used as a post-processing step for SOCP, akin to the previous analysis of running a DC OPF with an AC PF post-processing step. Note that other relaxations could be considered for comparison in future work; for example, the semidefinite programming relaxation (SDP). As \cite{QCRelax} mentions, the SDP relaxation is generally significantly slower than SOCP (in our simulations, the given computational platform was unable to solve SDP for the largest considered network), but provides a tighter relaxation. Tradeoffs are dependent on the application - the application here focuses on real time AC OPF.\\

\begin{table}[t!]\centering \small \caption{Percent error for the proposed method and for SOCP across the 500-sample testing sets.}
\begin{tabular}{|l|
>{\columncolor[HTML]{FFFFFF}}l |l|
>{\columncolor[HTML]{FFFFFF}}l |l|}
\hline
{\textbf{Network}}  & { \textbf{\begin{tabular}[c]{@{}l@{}}Average \\ Solve \\ Time \\ (NN)\vspace{2mm}\end{tabular}}} & { \textbf{\begin{tabular}[c]{@{}l@{}}Average \\ Solve\\ Time\\ (SOCP)\vspace{2mm}\end{tabular}}} & { \textbf{\begin{tabular}[c]{@{}l@{}}Average \\  Error \\ (NN)\end{tabular}}} & { \textbf{\begin{tabular}[c]{@{}l@{}}Average \\ Error \\ (SOCP)\end{tabular}}} \\ \hline
{ \textbf{30-bus}}   & {0.05 s}  & {0.15 s}      & {0.20\%}       & {2.86\%}      \\ \hline
{\textbf{500-bus}}  & {0.12 s}  & {3.48 s}    & {2.70\%}      & {26.85\%}      \\ \hline
{\textbf{1354-bus}} & {0.25 s}    & {33.72 s}            & {0.09\%}      & {0.28\%}                                                                        \\ \hline
\end{tabular}\label{tab:socp}
\end{table}

\section{Conclusion and Discussion} \label{sec:conclusion}
This paper provided a learning-based approximation for solving AC optimal power flow. Promising initial results indicate that the method can achieve very fast convergence speeds and high accuracy. Directions of future work include development of datasets or dataset generation methods for learning-based OPF, inclusion of additional constraints such as line limits, and speed/accuracy comparisons with other relaxations or convexifications. 

There are multiple drawbacks of this learning-based method that should be discussed here; although these results show an effective speed gain, traditional optimization has multiple upsides that are not yet covered by the ``learning-for-OPF'' models in literature (which also provide promising future directions of research). First, the model is trained on one static configuration of the network. That means that for any transformer tap-changing, any line switching, or switchable shunt control for example, a different model would have to be trained. Second, the performance of the model is limited by the dataset from which it was trained on. Lastly, neural networks do not offer grid operators as much insight or confidence in the resulting decision-making as traditional optimization does. This heuristic also does not provide a theoretical lower bound like other relaxations such as SOCP and SDP.

\section*{Acknowledgement and Data Availability}
This work is supported by the National Science Foundation CAREER award 2041835. The datasets and code used in this paper may be obtained by emailing the author.

% \noindent where $\bp_g$ is a vector comprising the active power generation $p_{g,j}$ at each generator $j \in \cal{G}$ and $\bv$ is $2n$-dimensional vector comprising the unknown voltage magnitudes and angles. It is clear from \eqref{eqn:PF_constr} and \eqref{eqn:PF_constr2} that the overall optimization problem is nonconvex. 

%%%%%%%%%%%%%%%%%%%%%%%%%%%%%%%%%%%%%%%%
%\bibliographystyle{IEEEbib}

\bibliographystyle{IEEEtran}
\bibliography{references.bib}

\begin{thebibliography}{10}
\providecommand{\url}[1]{#1}
\csname url@samestyle\endcsname
\providecommand{\newblock}{\relax}
\providecommand{\bibinfo}[2]{#2}
\providecommand{\BIBentrySTDinterwordspacing}{\spaceskip=0pt\relax}
\providecommand{\BIBentryALTinterwordstretchfactor}{4}
\providecommand{\BIBentryALTinterwordspacing}{\spaceskip=\fontdimen2\font plus
\BIBentryALTinterwordstretchfactor\fontdimen3\font minus
  \fontdimen4\font\relax}
\providecommand{\BIBforeignlanguage}[2]{{%
\expandafter\ifx\csname l@#1\endcsname\relax
\typeout{** WARNING: IEEEtran.bst: No hyphenation pattern has been}%
\typeout{** loaded for the language `#1'. Using the pattern for}%
\typeout{** the default language instead.}%
\else
\language=\csname l@#1\endcsname
\fi
#2}}
\providecommand{\BIBdecl}{\relax}
\BIBdecl

\bibitem{Zamzam_learn_19}
A.~Zamzam and K.~Baker, ``Learning optimal solutions for extremely fast {AC}
  optimal power flow,'' in \emph{IEEE SmartGridComm}, Dec. 2020.

\bibitem{chatzos2020highfidelity}
M.~Chatzos, F.~Fioretto, T.~W.~K. Mak, and P.~V. Hentenryck, ``High-fidelity
  machine learning approximations of large-scale optimal power flow,''
  \emph{arXiv preprint arXiv:2006.16356}, 2020.

\bibitem{DeepOPF_AC}
X.~Pan, M.~Chen, T.~Zhao, and S.~Low, ``Deep{OPF}: A feasibility- optimized
  deep neural network approach for {AC} optimal power flow problems,''
  \emph{arXiv preprint arXiv:2007.0100}, Jul. 2020.

\bibitem{MATPOWER}
R.~D. Zimmerman, C.~E. Murillo-Sanchez, and R.~J. Thomas, ``{MATPOWER}:
  Steady-state operations, planning, and analysis tools for power systems
  research and education,'' \emph{IEEE Trans. on Power Systems}, vol.~26,
  no.~1, pp. 12--19, Feb 2011.

\bibitem{PEGASE}
S.~{Fliscounakis}, P.~{Panciatici}, F.~{Capitanescu}, and L.~{Wehenkel},
  ``Contingency ranking with respect to overloads in very large power systems
  taking into account uncertainty, preventive, and corrective actions,''
  \emph{IEEE Trans. on Power Systems}, vol.~28, no.~4, pp. 4909--4917, 2013.

\bibitem{Baker_learning}
K.~Baker, ``Learning warm-start points for {AC} optimal power flow,'' in
  \emph{IEEE Machine Learning for Signal Proc. Conf. (MLSP)}, Oct. 2019.

\bibitem{chen2020hotstarting}
L.~Chen and J.~E. Tate, ``Hot-starting the {Ac} power flow with convolutional
  neural networks,'' \emph{arXiv preprint arXiv:2004.09342}, 2020.

\bibitem{homoSanja}
S.~{Cvijic}, P.~{Feldmann}, and M.~{Ilic}, ``Applications of homotopy for
  solving {AC} power flow and {AC} optimal power flow,'' in \emph{2012 IEEE
  Power and Energy Society General Meeting}, 2012, pp. 1--8.

\bibitem{powermodels}
C.~Coffrin, R.~Bent, K.~Sundar, Y.~Ng, and M.~Lubin, ``Powermodels.jl: An
  open-source framework for exploring power flow formulations,'' in \emph{Power
  Systems Computation Conference (PSCC)}, June 2018, pp. 1--8.

\bibitem{SOCPJabr}
R.~A. Jabr, ``Radial distribution load flow using conic programming,''
  \emph{IEEE Trans. on Power Sys.}, vol.~21, no.~3, pp. 1458--1459, Aug 2006.

\bibitem{QCRelax}
C.~Coffrin, H.~L. Hijazi, and P.~Van~Hentenryck, ``The qc relaxation: A
  theoretical and computational study on optimal power flow,'' \emph{IEEE
  Transactions on Power Systems}, vol.~31, no.~4, pp. 3008--3018, 2016.

\end{thebibliography}
%%%%%%%%%%%%%%%%%%%%%%%%%%%%%%%%%%%%%%%%

\end{document}